
\magnification=\magstep1
\def\to{\ \longrightarrow\ }

\def\nl{\hfill\break}

\def\hexnumber#1{\ifcase#1 0\or 1\or 2\or 3\or 4\or 5\or 6\or 7\or 8\or
 9\or A\or B\or C\or D\or E\or F\fi}
%
%
\font\twelvemsa=msam10 scaled 1200   
\font\tenmsa=msam10                  
\font\ninemsa=msam9            \font\sevenmsa=msam7
\font\sixmsa=msam6             \font\fivemsa=msam5
%
%
\newfam\msafam                 \textfont\msafam=\tenmsa
\scriptfont\msafam=\sevenmsa   \scriptscriptfont\msafam=\fivemsa
\edef\hexa{\hexnumber\msafam}        
\def\msa{\fam\msafam\tenmsa}         
%
%
\font\twelvemsb=msbm10 scaled 1200   
\font\tenmsb=msbm10                  
\font\ninemsb=msbm9            \font\sevenmsb=msbm7
\font\sixmsb=msbm6             \font\fivemsb=msbm5
%
\newfam\msbfam                 \textfont\msbfam=\tenmsb       
\scriptfont\msbfam=\sevenmsb   \scriptscriptfont\msbfam=\fivemsb
\edef\hexb{\hexnumber\msbfam}        
\def\msb{\fam\msbfam\tenmsb}         
%
%
\font\twelveeufm=eufm10 scaled 1200  
\font\teneufm=eufm10                 
\font\nineeufm=eufm9           \font\seveneufm=eufm7
\font\sixeufm=eufm6            \font\fiveeufm=eufm5
%
\newfam\eufmfam                \textfont\eufmfam=\teneufm
\scriptfont\eufmfam=\seveneufm \scriptscriptfont\eufmfam=\fiveeufm
\edef\hexf{\hexnumber\eufmfam}      
\def\frak{\fam\eufmfam\teneufm}     
%
%
%
\font\twelverm=cmr10 scaled 1200    
\font\ninerm=cmr9                   
\font\sixrm=cmr6   
%
\font\twelvei=cmmi10 scaled 1200    
\font\ninei=cmmi9                   
\font\sixi=cmmi6  
%
\font\twelvesy=cmsy10 scaled 1200   
\font\ninesy=cmsy9                  
\font\sixsy=cmsy6  
%
\font\twelvebf=cmbx10 scaled 1200   
\font\ninebf=cmbx9                  
\font\sixbf=cmbx6  
%
%
\font\twelveit=cmti10 scaled 1200   
\font\nineit=cmti9                  
%
\font\twelvesl=cmsl10 scaled 1200   
\font\ninesl=cmsl9                  
%
\font\twelvett=cmtt10 scaled 1200   
\font\ninett=cmtt9                  
%
%
%
%
\def\small{%
%
%
\textfont0=\ninerm \scriptfont0=\sixrm \scriptscriptfont0=\fiverm
\def\rm{\fam0\ninerm}        
%
%
\textfont1=\ninei \scriptfont1=\sixi \scriptscriptfont1=\fivei
%
%
\textfont2=\ninesy \scriptfont2=\sixsy \scriptscriptfont2=\fivesy
%
%
\textfont3=\tenex \scriptfont3=\tenex \scriptscriptfont3=\tenex
%
%
\textfont\bffam=\ninebf \scriptfont\bffam=\sixbf
\scriptscriptfont\bffam=\fivebf \def\bf{\fam\bffam\ninebf}%
%
%
\textfont\itfam=\nineit \def\it{\fam\itfam\nineit}%
\textfont\slfam=\ninesl \def\sl{\fam\slfam\ninesl}%
\textfont\ttfam=\ninett \def\tt{\fam\ttfam\ninett}%
%
%
%
\textfont\msafam=\ninemsa \scriptfont\msafam=\sixmsa
\scriptscriptfont\msafam=\fivemsa \def\msa{\fam\msafam\ninemsa}%
%
%
\textfont\msbfam=\ninemsb \scriptfont\msbfam=\sixmsb
\scriptscriptfont\msbfam=\fivemsb \def\msb{\fam\msbfam\ninemsb}%
%
%
\textfont\eufmfam=\nineeufm  \scriptfont\eufmfam=\sixeufm
\scriptscriptfont\eufmfam=\fiveeufm \def\frak{\fam\eufmfam\nineeufm}%
%
%
%
\normalbaselineskip=11pt
\setbox\strutbox=\hbox{\vrule height8pt depth3pt width0pt}%
%
%
\normalbaselines\rm}    
%
%
%
%
\def\large{%
\textfont0=\twelverm \scriptfont0=\ninerm \scriptscriptfont0=\sevenrm
\def\rm{\fam0\twelverm}%
\textfont1=\twelvei \scriptfont1=\ninei \scriptscriptfont1=\seveni
\textfont2=\twelvesy \scriptfont2=\ninesy \scriptscriptfont2=\sevensy
\textfont3=\tenex \scriptfont3=\tenex \scriptscriptfont3=\tenex
\textfont\bffam=\twelvebf \scriptfont\bffam=\ninebf
\scriptscriptfont\bffam=\sevenbf \def\bf{\fam\bffam\twelvebf}%
\textfont\itfam=\twelveit \def\it{\fam\itfam\twelveit}%
\textfont\slfam=\twelvesl \def\sl{\fam\slfam\twelvesl}%
\textfont\ttfam=\twelvett \def\tt{\fam\ttfam\twelvett}%
\textfont\msafam=\twelvemsa \scriptfont\msafam=\ninemsa
\scriptscriptfont\msafam=\sevenmsa \def\msa{\fam\msafam\twelvemsa}         
\textfont\msbfam=\twelvemsb \scriptfont\msbfam=\ninemsb
\scriptscriptfont\msbfam=\sevenmsb \def\msb{\fam\msbfam\twelvemsb}         
\textfont\eufmfam=\twelveeufm  \scriptfont\eufmfam=\nineeufm
\scriptscriptfont\eufmfam=\seveneufm \def\frak{\fam\eufmfam\teneufm}
\normalbaselineskip=15pt
\setbox\strutbox=\hbox{\vrule height11pt depth4pt width0pt}%
\normalbaselines\rm}%
%
\def\Bbb{\msb}

%

%
\mathchardef\plussquare="0\hexa01
\mathchardef\nge="3\hexb0B
\mathchardef\maltesecross="0\hexa7A
\mathchardef\del="0\hexf01
%

%

\input epsf
\overfullrule=0pt

\font\npt=cmr9
\font\Bbb=msbm10

\font\secfont=cmbx10
\font\ab=cmbx8

\font\nam=cmr8
\font\aff=cmti8

\font\em=cmti10

\mathchardef\square="0\hexa03
\def\qed{\hfill$\square$\par\rm}
\def\np{\vfill\eject}
\def\boxing#1{\ \lower 3.5pt\vbox{\vskip 3.5pt\hrule \hbox{\strut\vrule
\ #1 \vrule} \hrule} }

\def\down#1{\ \lower 3.5pt\vbox{\vskip 3.5pt \hbox{\strut \ #1 \vrule} \hrule} }
\def\negdown#1{\ \lower 3.5pt\vbox{\vskip 3.5pt \hbox{\strut  \vrule \ #1 }\hrule} }

\hsize=6.3 truein
\vsize=9 truein

\baselineskip=13 pt
\parskip=\baselineskip
 1

\parindent=0pt

\def\Z{\hbox{\Bbb Z}}
\def\R{\hbox{\Bbb R}}
\def\C{\hbox{\Bbb C}}
\def\H{\hbox{\Bbb H}}

\def\op{\buildrel o \over p}
\def\np{\buildrel n \over p}
\def\oM{\buildrel o \over {\cal M}}
\def\oP{\buildrel o \over {P}}

\def\nM{\buildrel n \over {\cal M}}
\def\nP{\buildrel n \over {P}}

\def\oover#1{\vbox{\ialign{##\crcr
{\npt o}\crcr\noalign{\kern 1pt\nointerlineskip}
$\hfil\displaystyle{#1}\hfil$\crcr}}}







\newif \iftitlepage \titlepagetrue

\def\diagram{\global\advance\diagramnumber by 1
$$\epsfbox{newlongfig.\number\diagramnumber}$$}
\def\ddiagram{\global\advance\diagramnumber by 1
\epsfbox{newlongfig.\number\diagramnumber}}

\newcount\diagramnumber
\diagramnumber=0

\newcount\secnum \secnum=0
\newcount\subsecnum
\newcount\defnum
\def\section#1{
                \vskip 10 pt
                \advance\secnum by 1 \subsecnum=0
                \leftline{\secfont \the\secnum \rm\quad#1}
                }

\def\subsection#1{
                \vskip 10 pt
                \advance\subsecnum by 1 
                \defnum=1
                \leftline{\secfont \the\secnum.\the\subsecnum\ \rm\quad #1}
                }

\def\definition{
                \advance\defnum by 1 
                \bf Definition 
\the\secnum .\the\defnum \rm \ 
                }

\def\lemma{
                \advance\defnum by 1 
                \par\bf Lemma  \the\secnum
.\the\defnum \rm \ \par
                }

\def\theorem{
                \advance\defnum by 1 
                \par\bf Theorem  \the\secnum
.\the\defnum \rm \ 
               }

\def\cite#1{
				\secfont [#1]
				\rm$\!\!\!$\nobreak
}

\vglue 20 pt
\centerline{\secfont New Invariants of Long Virtual Knots}
\medskip

\centerline{\nam Andrew Bartholomew${}^1$}
\centerline{\nam Roger Fenn${}^1$}
\centerline{\nam Naoko Kamada${}^2$}
\centerline{\nam Seiichi Kamada${}^3$}

\bigskip
\centerline{\nam ABSTRACT}
\leftskip=0.25 in
\rightskip=0.25in

{\ab This paper extends the construction of invariants for virtual knots to
virtual long knots and introduces two new invariant modules of virtual long
knots.  Several interesting features are described that distinguish virtual
long knots from their classical counterparts with respect to their
symmetries and the concatenation product.}

\leftskip=0 in
\rightskip=0in

\centerline{\aff ${}^1$School of Mathematical Sciences, University of Sussex}
\centerline{\aff Falmer, Brighton, BN1 9RH, England}
\centerline{\aff e-mail address: rogerf@sussex.ac.uk}
\centerline{\aff ${}^2$This research is supported by the 21st COE program}
\centerline{\aff ``Constitution of wide-angle mathematical basis focused on knots''}
\centerline{\aff Department of Mathematics, Osaka City University, Osaka 558-8585, Japan}
\centerline{\aff e-mail address: naoko@sci.osaka-cu.ac.jp} 
\centerline{\aff ${}^3$This research is partially supported by Grant-in-Aid
for Scientific Research, JSPS}
\centerline{\aff  Department of Mathematics, Hiroshima University, Hiroshima 739-8526, Japan)} 
 \centerline{\aff e-mail address:  kamada@math.sci.hiroshima-u.ac.jp} 

\section{Introduction}

Virtual knots are a generalization of classical knots introduced by
L. Kauffman in 1996 \cite{kauD}.  They describe all knots in
thickened closed orientable surfaces of any genera, and are in one-to-one correspondence
with abstract equivalence classes (or stable equivalence classes)
\cite{kk, kauD, CKS2002, Kuper}.  They are a supplement to real
knots so that all Gauss codes (or Gauss diagrams) can be realized in
the category of virtual knots.  Consequently they are helpful for the
study of some invariants, including the Jones polynomials \cite{kamD,
kamDD, kamE, kk, kauD} and finite type invariants, \cite{GPV}.
They can be also used in rack and quandle homology theory when 
describing a $2$-cycle as a diagram \cite{CJKS2001a, CKS2001, FRSa,
FRSb, FRSc, FRSd, FRSe, Greene}.  Some invariants of classical
knots can be generalized to those of virtual knots, and the others
cannot.  Virtual knots also have their own invariants, like the JKSS
polynomials (Jaeger, Kauffman, Saleur \cite{JKS} and Sawollek
\cite{Saw}), Silver-Williams invariants \cite{SWA}, etc.  and the
quaternionic invariants \cite{BaF, BuF}.

A ({\it long}) {\it virtual knot diagram} is an oriented (long) knot
diagram which may have encircled crossings called virtual
crossings. Two diagrams are said to be {\it equivalent} if they are
related by a finite sequence of generalized Reidemeister moves
introduced in \cite{kauD}.  A ({\it long}) {\it virtual knot} is the
equivalence class of a (long) virtual knot diagram.  By closing the
ends of a long virtual knot diagram, we obtain a virtual knot diagram.
It induces a map from the set of long virtual knots to the set of
virtual knots.  However, unlike in the classical case, this map is not
a bijection or even injective (see \cite{GPV}).

In this paper we introduce four modules ${\cal M}_K$, $\widehat{\cal
M}_K$, ${\oM}_K$ and ${\nM}_K$ which are invariants of long knots. The
modules are left ${\cal F}$ modules where ${\cal F}$ is the algebra
introduced in \cite{BuF, BaF, F}. This has two generators $A, B$ and one relation
$$A^{-1}B^{-1}AB-B^{-1}AB=BA^{-1}B^{-1}A-A.$$
As in these previous papers, this algebra
can be represented onto more tractable algebras, say the quaternions,
and invariants, usually polynomials, can be calculated from these.
We do this by defining codimension $r$ determinants $\det_\eta^{(r)}(P)$ of a presentation matrix $P$ with representation $\eta$.

The definition of the module ${\cal M}_K$ just mimics that of the
module of a virtual knot, and $\widehat{\cal M}_K$ is the module of
the virtual knot obtained from $K$ by closing the ends, see \cite{F}.
The module ${\oM}_K$ is obtained by putting the input generator
equal to zero. The module ${\nM}_K$ is obtained by putting the output generator equal to zero. Thus ${\oM}_K$ and ${\nM}_K$  are the definitions that are
essentially new in this paper and have an interesting feature.  In
fact, the determinants of these modules satisfies a
product formula with respect to the concatenation product of two long
virtual knots.

The second named author, R. Fenn introduced the Budapest switch
(\cite{BuF, rFJK}), augmented by $t$, that is the $2 \times 2$ matrix
$S$ defined by
$$
S = \pmatrix{
1+i  &  -t j \cr
t^{-1} j  & 1+i \cr
}$$
where $i, j$ have the usual meanings as quaternions and $t$ is a
central variable.  The matrix is invertible and satisfies the set
theoretic Yang-Baxter equation (cf. \cite{BaF, BuF}).  It follows
that the first row entries, $1+i$ and $-tj$, define a representation
of the algebra ${\cal F}$.

In the cited papers, a method is described which defines a presentation
matrix derived from a diagram of a virtual knot $K$ and this
determines an invariant {\it quaternionic module} of $K$.  The
codimension $r$ Study determinants of the presentation matrix are 
invariants.  The following examples are from \cite{BaF, BuF}. For the virtual trefoil (the first figure in p.~24 of \cite{BuF}), the quaternionic module has a presentation matrix
$$
\pmatrix{
- t^2 + 2i  &  
-1 + t (-j +k) + t^{-1} (j + k) \cr
-1 + t ( -j -k) + t^{-1} (j-k)  & 
-t^{-2} + 2i\cr
}$$
(There is a typo in the (2,2)-entry of the matrix in p.~24 of
\cite{BuF}.)  The determinant is $1 + 2 t^2 + t^4$.
For the Kishino knot (the second figure in p.~24 of \cite{BuF}),
the determinant is  $0$ and the codimension 1 determinant is $1 + (5/2) t^2 + t^4$.

In this paper we shall use this and other representations to
define invariant polynomials from the four modules.

We can also repeat this analysis for {\em flat} long knots. These are
represented like long virtual knots with virtual crossings but instead of standard crossings they have flat crossings which are the projections of standard
crossings. These may be interpreted as paths on an oriented
surface. Once again we have four modules but these are now left
modules over the Weyl algebra, \cite{FT}.

The second named author would like to thank the third and the fourth
named authors and the Universities of Osaka and Kyoto for their
hospitality during his stay in March, 2006.

The third and the fourth named authors would like to thank the second named author
and University of Sussex for his hospitality during their stay in
September, 2005.

\section{Determinants from non-commuting rings}

In this section we show how determinants of matrices with entries in a
general ring $R$ can be defined. Here $R$ is a (possibly
non-commutative), associative ring.  The definition will depend upon a
representation of $R$ into the ring,  $M_{d,d}\Lambda$ of $d\times d$ square matrices with entries in a commutative ring $\Lambda$. For positive codimension determinants it will be useful  if $\Lambda$ supports a greatest common divisor function, written gcd. A good reference for this is \cite{As}.

As an illustration, consider the
non-commutative ring $R = {\H}[t, t^{-1}]$ of Laurent polynomials
whose coefficients are quaternions and the varible $t$ is central. The commutative ring is $\Lambda={\C}[t, t^{-1}]$, the Laurent polynomials in $t$ with complex coeficients. The representation
$\mu: R \to {M}_{2,2}( {\C}[t, t^{-1}])$ has $d=2$ and is
defined by
$$(\alpha_1  + \alpha_2 i + \alpha_3 j + \alpha_4 k) t^s \mapsto 
\pmatrix{(\alpha_1 + \alpha_2 i)  t^s &  (\alpha_3 + \alpha_4 i)   t^s  \cr
(-\alpha_3 + \alpha_4 i)  t^s & (\alpha_1 -\alpha_2 i)  t^s  \cr}
$$
where $\alpha_1, \dots, \alpha_4 \in {\R}$ and $s  \in {\Z}$. 
This representation is {\bf standard} and is usually used in illustrations. 

The ring ${\C}[t, t^{-1}]$ has greatest common divisors, see
\cite{CF}.

Returning now to the general case, let $\eta: R\to M_{d,d}\Lambda$ be
the representation and let $P\in M_{n,m}(R)$. A square submatrix $B$
of $P$ is said to have {\rm codimension} $r$ if it is obtained by
deleting $n-m+r$ rows and $r$ columns if $n\ge m$ or by deleting
$m-n+r$ columns and $r$ rows if $m\ge n$. For simplicity assume $m\ge n$.  Let $B_1, \dots, B_s$ be the codimension $r$ submatrices of $P$ of size $(n-r)\times (n-r)$. Consider $\eta(B_1), \dots, \eta(B_s)$,
which are $d(n-r)\times d(n-r)$ matrices whose entries belong to
$\Lambda$.  The {\it codimension $r$ $\eta$-determinant}, of $P$ is
the gcd of the usual determinants of these matrices.  We denote it by
$\det_\eta^{(r)}(P)$. If $P$ is square then $\det_\eta^{(0)}(P)$ is defined even if $\Lambda$ does not possess greatest common divisors. All determinants are well defined up to multiplication by a
unit.

Now suppose that ${\cal M}$ is a finitely presented $R$-module with
presentation matrix $P\in M_{m,n}(R)$. Then the elements
$\det_\eta^{(r)}(P)$ will be invariants of the module.

\section{The Invariant Modules}

Suppose that $R$ is an associative ring
and
$$S=\pmatrix{A&B\cr C&D\cr}$$
is a $2 \times 2$ matrix with entries from $R$. If $S$ is invertible
and satisfies the set theoretic Yang-Baxter equation in the sense of
\cite{BaF, FJK} then $S$ is called a {\em switch}. The
universal case occurs when $R={\cal F}$ and ${\cal F}$ has the
presentation
$${\cal F}=<A,B\mid A^{-1}B^{-1}AB-B^{-1}AB=BA^{-1}B^{-1}A-A>$$
and $C, D$ are defined by
$$C=A^{-1}B^{-1}A(1-A),\quad D=1-A^{-1}B^{-1}AB.$$

Let $x_0, x_1, \dots, x_n$ be the semi-arcs of a diagram of a long
virtual knot $K$, which appear in this order along $K$.  For each
positive crossing, we consider a relation
$$
S 
\pmatrix{
x_i \cr
x_j\cr}
=
\pmatrix{x_{j+1} \cr
x_{i+1}\cr}
$$
where $x_i$ and $x_j$ are incoming semi-arcs and $x_{j+1}$ and
$x_{i+1}$ are outgoing semi-arcs such that $x_j$ and $x_{j+1}$ are
under-arcs.  For each negative crossing, we consider a relation that
is the inverse of the positive one . 
$$
S 
\pmatrix{
x_{j+1} \cr
x_{i+1}\cr}
=
\pmatrix{x_{i} \cr
x_{j}\cr}
$$
Of course for a virtual crossing the labeling carries over. See the following diagram, (cf. \cite{BaF, FJK}).
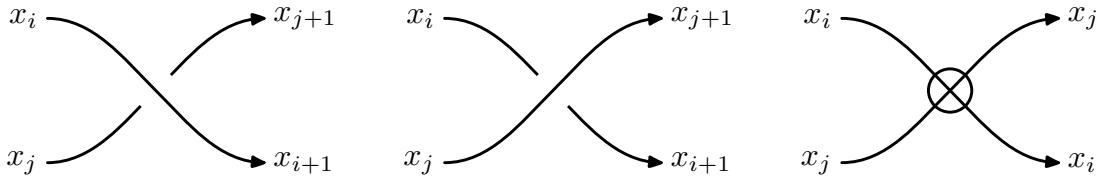
\diagram
The corresponding relations are
$$x_{j+1}=Ax_i+Bx_j,\quad x_{i+1}=Cx_i+Dx_j$$
and
$$x_{i}=Ax_{j+1}+Bx_{i+1},\quad x_{j}=Cx_{j+1}+Dx_{i+1}$$
The module ${\cal M}_K$ is the $R$-module generated by $x_0, x_1,
\dots, x_n$ with the relations associated with positive crossings and
negative crossings.  There is one more generator than relation.

The module $\widehat{\cal M}_K$ is the quotient of ${\cal M}_K$ by an
additional relation $x_0=x_n$. The number of generators and relations
are the same. That is, the presentation matrix is square.

The module ${\oM}_K$ is the quotient of ${\cal M}_K$ by an
additional relation $x_0=0$. Again, the presentation matrix is square.

The module ${\nM}_K$ is the quotient of ${\cal M}_K$ by an
additional relation $x_n=0$. The generator $x_n$ is the label on the outgoing arc. Again, the presentation matrix is square.

\theorem{
Suppose $A, 1-A$ and $B$ are invertible. Then the modules ${\cal M}_K$,
$\widehat{\cal M}_K$, ${\oM}_K$ and  ${\nM}_K$ are invariants of a long virtual knot $K$.}

{\it Proof}: Since $S$ is invertible, satisfies the set theoretic
Yang-Baxter equation, and since $1-A$ is invertible these module are
preserved by all generalized Reidemeister moves (see \cite{BaF,
FJK}).  \qed

As an illustration consider the ``fly'' long knot, $F$,  pictured below.

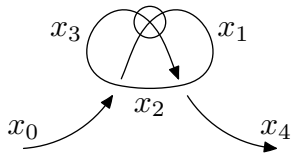
\diagram

The presentation matrices in the four cases are
$${\cal M}=\pmatrix{
-1&B&0&A&0\cr
0&D&-1&C&0\cr
0&A&0&B&-1\cr
0&C&-1&D&0\cr},\
\widehat{\cal M}=\pmatrix{
-1&B&0&A&0\cr
0&D&-1&C&0\cr
0&A&0&B&-1\cr
0&C&-1&D&0\cr
1&0&0&0&-1\cr}$$
and
$$
{\oM}=\pmatrix{
-1&B&0&A&0\cr
0&D&-1&C&0\cr
0&A&0&B&-1\cr
0&C&-1&D&0\cr
1&0&0&0&0\cr},\
{\nM}=\pmatrix{
-1&B&0&A&0\cr
0&D&-1&C&0\cr
0&A&0&B&-1\cr
0&C&-1&D&0\cr
0&0&0&0&1\cr},\
$$
\section{The Invariant Polynomials}

Let $K$ be a long knot and let ${\cal M}_K$, $\widehat{\cal M}_K$,
${\oM}_K$  and ${\nM}_K$ be the ${\cal F}$-modules defined in the previous section. Let $P$,
$\widehat P$,  $\oP$ and $\nP$ be the respective presentation
matrices. Suppose we now represent the algebra as matrices so that we can define determinantal invariants as described in section 2. Each entry in the $P$
matrices is a $d\times d$ matrix for some $d$. Let $p^{(r)}_K$,
$\widehat{p}^{(r)}_K$, 
$\op^{\raise -8pt\hbox{$\scriptstyle(r)$}}_K$ 
and 
$\np^{\raise -8pt\hbox{$\scriptstyle(r)$}}_K$ 
be the corresponding
determinants in the commutative ring $\Lambda$, with codimension $r=0,1,2\ldots$.

Let $\widehat{K}$ be the closure of the long knot $K$. Then
$\widehat{K}$ also has an invariant ${\cal F}$-module, see \cite{F}.
Let $q^{(r)}_{\widehat{K}}$ be the sequence of determinants in $\Lambda$,
corresponding to the presentation.

\theorem{
\parindent=20pt
\item{ (1)} $q^{(r)}_{\widehat{K}}
=\widehat{p}^{(r)}_K$ 
\item{ (2)} $p^{(r)}_K$ divides $\op^{\raise -8pt\hbox{$\scriptstyle(r)$}}_K$. 
\item{ (3)} $p^{(r)}_K$ divides $\np^{\raise -8pt\hbox{$\scriptstyle(r)$}}_K$. 
}
\parindent=0pt

{\it Proof}:
\def\oQ{\buildrel o \over {Q}}

 
(1) This follows since $\widehat{\cal M}_K$ is
equal to the module of
the closure $\widehat{K}$ of $K$.
 
(2) The module ${\cal M}_K$ has an $n \times (n+1)$
matrix $P$ as a
presentation matrix such that the first column corresponds to $x_0$
and the last column corresponds to $x_n$. A codimension $r$ submatrix
is obtained from $P$ by deleting $r$ rows and $r+1$ columns. Let $B_1,
\dots, B_s$ be the codimension $r$ submatrices of $P$. Then
$p^{(r)}_K$ divides the determinant of all of these after the
representation and is the largest, by division, element which does so.
 
The module ${\oM}_K$ has an $(n+1) \times (n+1)$
presentation matrix $\oP$ that is obtained from $P$ by adding
the row $(1, 0, \dots, 0)$ to the bottom. 
The presentation matrix $\oP$ is simplified to $\oQ$ that is an
$n \times n$ matrix obtained from $\oP$ by deleting the first column  and the bottom row, which  is obtained from $P$ by deleting the first column.
Since $\oQ$ is
square a codimension $r$ submatrix is obtained from $\oQ$ by
deleting $r$ rows and $r$ columns. It follows that
$\op^{\raise -8pt\hbox{$\scriptstyle(r)$}}_K$
 is the gcd of a subset of the values for which
$p^{(r)}_K$ is the gcd. The result now follows.
 
(3) The proof is similar to (2). \qed
 

\theorem{ Let $K_1 \cdot K_2$ be the concatenation product of
two long virtual knots $K_1$ and $K_2$. For any suitable
representation of the fundamental modules
$$\op^{\raise -8pt\hbox{$\scriptstyle(0)$}}(K_1 \cdot K_2) = 
\op^{\raise -8pt\hbox{$\scriptstyle(0)$}}(K_1) 
\op^{\raise -8pt\hbox{$\scriptstyle(0)$}}(K_2)$$
and
$$\np^{\raise -8pt\hbox{$\scriptstyle(0)$}}(K_1 \cdot K_2) = 
\np^{\raise -8pt\hbox{$\scriptstyle(0)$}}(K_1) 
\np^{\raise -8pt\hbox{$\scriptstyle(0)$}}(K_2). $$
}

{\it Proof}: Let $P_1= (a_0 a_1 \cdots a_n)$ and $P_2 =(b_0 b_1 \cdots
b_{n'})$ be the presentation matrices of ${\cal M}_{K_1}$ and ${\cal
M}_{K_2}$ associated with their diagrams.  Then ${\oM}_{K_1}$,
${\oM}_{K_2}$ and ${\oM}_{K_1 \cdot K_2}$ have presentation matrices
$$
 (a_1  \cdots  a_n), \quad  (b_1 \cdots  b_{n'}), \quad 
\pmatrix{
a_1  \cdots  a_{n-1}  &      a_n        & 0 \cdots 0 \cr
  0 \cdots   0              &           b_0  &   b_1  \cdots   b_{n'} \cr}
$$
respectively.  Thus we have the result for $\op$. The proof for $\np$ is
similar.
\qed
\section{Simplifying the Modules and some Calculations}
The presentation matrices defined above can be simplified by the usual rules for manipulating non-commuting relations. That is
{\parindent=20pt
\item{1} Interchange any row(column).
\item{2} Multiply any row(column) on the left(right) by a unit.
\item{3} Add any row(column) multiplied on the left(right) to a different row(column).
\item{4} Introduce or delete any zero row.
\item{5} $P\leftrightarrow \pmatrix{1&0\cr 0& P\cr}$
}

If we apply these rules to the example in section 3 and simplify as far as possible we get,
$$\pmatrix{
D-C&D-C\cr},\
\pmatrix{0\cr},\
\pmatrix{(C-D)(1+B^{-1}A)\cr},\
\pmatrix{(D-C)(1+B^{-1}A)\cr}$$
for the four presentation matrices. Note that the presentation matrix for $\widehat{\cal M}_K$ will reduce to zero since the closure of $B$ is the trivial knot.

Now apply the homomorphism which replaces $A, B, C, D$ with $1+i, -tj, t^{-1} j, 1+i$ respectively and use the standard representation. The three codimension zero polynomial invariants are
$|D-C|^2=2+t^{-2}$, $0$ and $|(C-D)(1+B^{-1}A)|^2=(2+t^{-2})(1+2t^{-2})$. All the higher codimension polynomials are 1.
\section{Symmetries of Long Virtual Knots}
There are various symmetries of the knot diagram which can be applied.
Consider reflection in the plane of a knot diagram $D$. Let $-D$ denote the resulting diagram. This interchanges plus and minus crossings. Let $\overline D$ denote the result of reflection in the $x$-axis. Finally let $D^*$ be obtained by reversing the arrow and rotating the result through 180 degrees. 

For the fly $F^*=\overline F$.
The effect of the other three possibilities on the fly are illustrated below.
\medskip
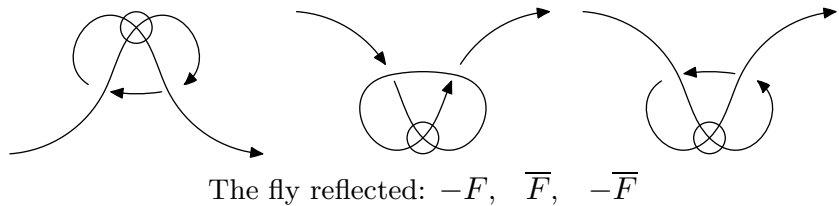
\centerline{\ddiagram\quad \ddiagram\quad \ddiagram}
\centerline{\npt The fly reflected: $-F$,\quad $\overline F$,\quad $-\overline F$}
Using a suitable representation of the fundamental algebra all three can be distinguished from themselves and the original fly as follows.

The resulting polynomials are tabulated as follows. The switch used is given by a representation of the quantum Weyl
algebra, with  $A, B, C, D$ the following $2\times2$ matrices.

$$
	A = \pmatrix{1-q & -q^3+2q^2-1 \cr 0 & 1-q },\ 
	B = \pmatrix{q & 1 \cr 0 & 1} 
$$
$$
	C = \pmatrix{1 & (-q^4+3q^3-2q^2-2q+1)/q \cr 0 & q},\ 
	D = \pmatrix{0 & (q^3-2q^2+1)/q \cr 0 & 0 } 
$$

$$\vbox{ \offinterlineskip \halign{
     \vrule # & \hfil \quad $#$ \quad \hfil & 
     \vrule # & \hfil     \ $#$ \     \hfil & 
     \vrule # & \hfil     \ $#$ \     \hfil & 
     \vrule # & \hfil     \ $#$ \     \hfil &
     \vrule # & \hfil     \ $#$ \     \hfil &
	 \vrule#\cr 
\noalign{\hrule} 
& K && p^{(0)}(K) && \widehat{p}^{(0)}(K) && \op^{\raise
-8pt\hbox{$\scriptstyle(0)$}}(K) && \np^{\raise
-8pt\hbox{$\scriptstyle(0)$}}(K) & \cr 
\noalign {\hrule} 
height1pt&\omit&height1pt&\omit&height1pt&\omit&height1pt&\omit&height1pt&\omit&height1pt\cr
& F           && 1 && 0 && (2-q)/q && (2-q)/q & \cr 
\noalign {\hrule} 
height1pt&\omit&height1pt&\omit&height1pt&\omit&height1pt&\omit&height1pt&\omit&height1pt\cr
&-F           && 2-q && 0 && (2-q)/q && (2-q)/q & \cr 
\noalign {\hrule} 
height1pt&\omit&height1pt&\omit&height1pt&\omit&height1pt&\omit&height1pt&\omit&height1pt\cr
&\overline F  && 1  && 0 && 2q-1 && 2q-1 & \cr
\noalign {\hrule} 
height1pt&\omit&height1pt&\omit&height1pt&\omit&height1pt&\omit&height1pt&\omit&height1pt\cr
&-\overline F &&  2q-1 && 0 && 2q-1 && 2q-1 &\cr 
\noalign {\hrule} 
}}$$


\theorem{Consider the following three conditions on a switch $S$.
a) $S=S^{\dag}$, b) $S^2=1$, c) $SS^{\dag}=1$, where
$S^{\dag}=\pmatrix{D&C\cr C&A\cr}$.
\nl
If a) then $A=D$ and $B=C$ and $K$ cannot be distinguished from $-\overline K$. \nl
If b) then the underlying algebra is the Weyl algebra and $K$ cannot be distinguished from $-K$ or $K^*$. \nl
If c) then $A, B$ commute and $S$ is $\pmatrix{2&\pm1\cr \mp1&0\cr}$, a specialization of the Alexander switch. Moreover $K$ cannot be distinguished from $\overline K$.}

{\it Proof}: Most of the results easily follow by looking at the conditions on the entries of $S$ and how this affects the calculations of the modules.
\nl
For b) the underlying algebra is the Weyl algebra because of results in \cite{FT}.
\nl
For c) the condition $SS^{\dag}=1$ implies 
$$\matrix{
AD+B^2&=\hfill1&AC&=-BA \cr CD&=-DB&C^2+DA&=\hfill1\cr}.$$
Using the relations
$$C=A^{-1}B^{-1}A(1-A),\quad D=1-A^{-1}B^{-1}AB$$
gives the result. \qed 

It is well known that the product of two classical knots is a commutative operation. This is not the case for the product of two long virtual knots.
For example the products $F\cdot \overline F$ and $\overline F\cdot F$ are distinct. 
A calculation using the Budapest switch shows that for $F\cdot \overline F$, $ 
p^{(0)} = 6t^4+15t^2+6$ whereas for $\overline F\cdot F$ we have
$p^{(0)} = 3t^4+15/2t^2+3$, which is half the previous polynomial. If we are working over the integer quaternions then 2 is not a unit and so this shows the knots are distinct. 

This is not perhaps a "killer" example. If we consider
$$(F \cdot \overline{F}) \cdot F \hbox{ and  }F\cdot (F \cdot \overline{F})$$
then in both cases $p^{(0)} =-12t^8-60t^6-99t^4-60t^2-12$ but the values of $p^{(1)} $ are $9(t^2+1)$ and $-9/2(2t^2+1)(t^2+2)$ respectively.

\section{Long flat virtual knots}  

We now repeat the previous analysis for long flat knots. For a full
discusion of the following see \cite{FT}. A switch $S$ can be used
provided it satisfies $S^2=id$. The conditions for this are contained
in the following theorem.

\theorem{
Suppose $A, 1-A$ and $B$ are invertible and
$$
S = 
\pmatrix{
A& B \cr
C& D\cr
}$$
is a $2 \times 2$ matrix with entries satisfying
$$A^{-1}B^{-1}AB-B^{-1}AB=BA^{-1}B^{-1}A-A=1$$
and $C, D$ are defined by
$$C=A^{-1}B^{-1}A(1-A),\quad D=1-A^{-1}B^{-1}AB.$$
Then $S^2=1$ if and only if $u=B, v=B^{-1}A^{-1}$ satisfy $uv-vu=1$.
(The elements $u, v$ are the generators of the Weyl algebra.)}
\qed

We now repeat the construction considered earlier and arrive at modules ${\cal WM}_F$, $\widehat{\cal WM}_F$, ${\cal W}\!\oM_F$ and ${\cal W}\!\nM_F$ which are invariants of a long flat knot $F$.

Again by analogy we can find tractible invariants given representations
onto finite matrices. Many examples are given in \cite{FT}.

The following is an example with ring ${\Z}_2[a, x, y]$.
$$
u = \pmatrix{
x  &  a \cr
0  &  x \cr},\
v = \pmatrix{
y  &  0 \cr
1/a  &  y \cr}
$$
Consider the "flat fly" denoted by $FF$ and its reflection in the $x$-axis illustrated below.
\medskip
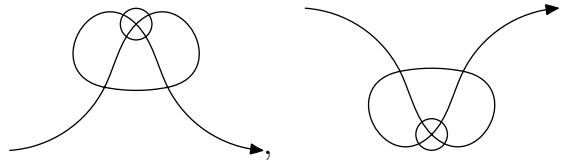
\centerline{\ddiagram,\quad\ddiagram}

\centerline{$FF$ and $\overline{FF}$}

Using the representation above, setting $y=x$ and $a=1$, the codimension zero polynomials are

$$\vbox{ \offinterlineskip \halign{
     \vrule # & \hfil \quad $#$ \quad \hfil & 
     \vrule # & \hfil     \ $#$ \     \hfil & 
     \vrule # & \hfil     \ $#$ \     \hfil & 
     \vrule # & \hfil     \ $#$ \     \hfil &
     \vrule # & \hfil     \ $#$ \     \hfil &
	 \vrule#\cr 
\noalign{\hrule} 
& K && p^{(0)}(K) && \widehat{p}^{(0)}(K) && \op^{\raise
-8pt\hbox{$\scriptstyle(0)$}}(K) && \np^{\raise
-8pt\hbox{$\scriptstyle(0)$}}(K) & \cr 
\noalign {\hrule} 
height1pt&\omit&height1pt&\omit&height1pt&\omit&height1pt&\omit&height1pt&\omit&height1pt\cr
& {\bf d}(FF)  && x^2+1 && 0 && x^8+x^6+x^2+1 && x^8+x^6+x^2+1 & \cr 
\noalign {\hrule} 
height1pt&\omit&height1pt&\omit&height1pt&\omit&height1pt&\omit&height1pt&\omit&height1pt\cr
&{\bf d}(\overline{FF}) && x^6+1 && 0 && x^8+x^6+x^2+1 && x^8+x^6+x^2+1 & \cr 
\noalign {\hrule} 
}}$$

This shows that $FF$ and $\overline{FF}$ are both non-trivial and distinct.

It is interesting to note Turaev's descent map ${\bf d}$ of
long flat knots to long virtual knots as an alternative method to show that $FF$ is non-trivial. This lifts flat knots by turning the first time a crossing is met to an overcrossing. For example $FF$ and $\overline{FF}$ are converted as shown in the following diagram.

$$\matrix{\ddiagram\ &\raise 20 pt\hbox{$\to\ $} &\ddiagram\cr}
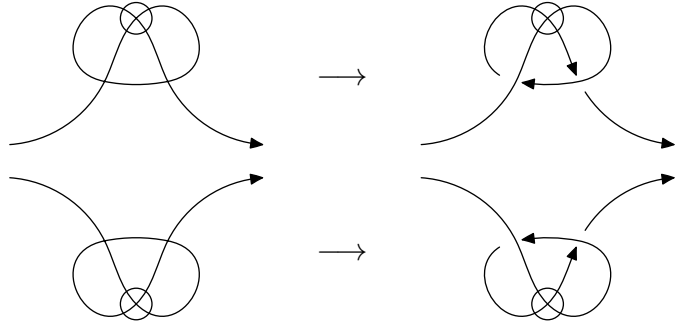
$$
$$\matrix{\ddiagram\ &\raise 20 pt\hbox{$\to\ $} &\ddiagram\cr}$$

Then, using the same switch as in the previous example, the polynomials for ${\bf d}(FF)$ and ${\bf d}(\overline{FF})$ 
are also

$$\vbox{ \offinterlineskip \halign{
     \vrule # & \hfil \quad $#$ \quad \hfil & 
     \vrule # & \hfil     \ $#$ \     \hfil & 
     \vrule # & \hfil     \ $#$ \     \hfil & 
     \vrule # & \hfil     \ $#$ \     \hfil &
     \vrule # & \hfil     \ $#$ \     \hfil &
	 \vrule#\cr 
\noalign{\hrule} 
& K && p^{(0)}(K) && \widehat{p}^{(0)}(K) && \op^{\raise
-8pt\hbox{$\scriptstyle(0)$}}(K) && \np^{\raise
-8pt\hbox{$\scriptstyle(0)$}}(K) & \cr 
\noalign {\hrule} 
height1pt&\omit&height1pt&\omit&height1pt&\omit&height1pt&\omit&height1pt&\omit&height1pt\cr
& {\bf d}(FF)  && x^2+1 && 0 && x^8+x^6+x^2+1 && x^8+x^6+x^2+1 & \cr 
\noalign {\hrule} 
height1pt&\omit&height1pt&\omit&height1pt&\omit&height1pt&\omit&height1pt&\omit&height1pt\cr
&{\bf d}(\overline{FF}) && x^6+1 && 0 && x^8+x^6+x^2+1 && x^8+x^6+x^2+1 & \cr 
\noalign {\hrule} 
}}$$

Our invariants are consequence of biquandles \cite{FJK}.  When we deform a given virtual knot or a given long virtual knot into a braid form, it is easier to calculate the biquandle, the quaternionic module and the determinant invariants.  For braiding of virtual knots and long virtual knots, refer to 
\cite{SkamA, SkamB, kaulam}.  
\section{References}

{\bf As} H. Aslaksen: 
{\it Quaternionic determinants}, 
Math. Intel. {\bf 18} (1996), 1-19

{\bf BaF}
A. Bartholomew and R. Fenn: 
{\it Quaternionic invariants of virtual knots and links}, 
preprint. 

{\bf BuF}
S. Budden and R. Fenn: 
{\it The equation $[B, (A-1)(A,B)]=0$ and virtual knots and links}, 
Fund. Math. {\bf 184} (2004), 19 --29. 


{\bf CJKS2001a} 
J. S. Carter, D. Jelsovsky, S. Kamada and M. Saito: 
  {\it Quandle homology groups, their Betti numbers, and virtual knots}, 
J. Pure Appl. Algebra {\bf 157} (2001), 135--155. 


{\bf CKS2001} 
J. S. Carter, S. Kamada and M.  Saito:  
{\it Geometric interpretations of quandle homology},  
J. Knot Theory Ramifications {\bf 10} (2001), 345--386. 

{\bf CKS2002} 
J. S. Carter, S. Kamada and M.  Saito:  
{\it  Stable equivalence of knots on surfaces and virtual knot cobordisms},  
J. Knot Theory Ramifications {\bf 11} (2002), 311--322. 

{\bf CF}
R. H. Crowell and R. H. Fox: 
{\it An introduction to knot theory}, 
Ginn and Co, (1963).

{\bf F}
R. Fenn: 
{\it Quaternion algebras and invariants of virtual knots and links I}, 
to appear in JKTR

{\bf FJK}
R. Fenn, M. Jordan, and L. Kauffman: 
{\it Biquandles and virtual links}, 
Topology Appl. {\bf 145} (2004), 157--175.

{\bf FT}
R. Fenn and V. Turaev: 
{\it Weyl Algebras and Knots}, 
J. Geometry and Physics {\bf 57} (2007), 1313-1324


{\bf FRSa}
R. Fenn, C. Rourke and B. Sanderson: 
{\it An introduction to species and the rack space}, 
``Topics in Knot Theory'' (M. E. Bozhuyu, ed.), Kluwer Academic, pp. 33--55, (1993). 

 {\bf FRSb}
R. Fenn, C. Rourke and B. Sanderson: 
{\it Trunks and classifying spaces,}
Appl. Categ. Structures {\bf 3} (1995), 321--356.

{\bf FRSc} 
R. Fenn, C. Rourke and B. Sanderson: 
{\it James bundles and applications,} preprint (1996), 
available
at
 http://www.maths.warwick.ac.uk/~cpr/ftp/james.ps

{\bf FRSd} 
R. Fenn, C. Rourke and B. Sanderson: 
{\it The rack space,} TAMS, 359 (2007) 701-740\nl
preprint (2003), 
available
at
arXiv: math.GT/0304228

{\bf FRSe}
R. Fenn, C. Rourke and B. Sanderson: 
{\it James bundles}, 
Proc. London Math. Soc. {\bf 89} (2004), 217--240.

{\bf GPV}
M. Goussarov, M. Polyak, and O. Viro: 
{\it Finite-type invariants of classical and virtual knots\/}, 
Topology {\bf 39} (2000), 1045--1068.

{\bf Greene} 
M. T. Greene: 
{\it Some results in geometric topology and geometry,}
Ph.D. Dissertation, Warwick (1997).


{\bf JKS}
F. Jaeger, L. H. Kauffman, and H. Saleur: 
{\it The conway polynomial in $R^3$ and in thickened surfaces: 
A new determinant formulation}, 
J. Combin. Theory Ser. B 
{\bf 61} (1994),  237--259.  





{\bf kamD} N. Kamada:  
{\it The crossing number of alternating link diagrams of a surface\/}, 
Proceedings of Knots 96, World Scientific Publishing Co.,1997, 377--382.


{\bf kamDD}
N.~Kamada:
{\it Span of the Jones polynomial of an alternating virtual link \/},
Algebr. Geom. Topol. {\bf 4} (2004), 1083--1101.


{\bf kamE} 
N.~Kamada: 
{\it A relation of Kauffman's $f$-polynomials of 
virtual links\/}, 
Topology and its Application {\bf 146--147} (2005), 123--132.

{\bf kk} 
N.~Kamada and S.~Kamada: 
{\it Abstract link diagrams and virtual knots\/}, 
J. Knot Theory Ramifications {\bf 9} (2000),
93--106.


 
{\bf SkamA}
S. Kamada:
{\it Braid presentation of virtual knots and welded knots},
Osaka J. Math., to appear,
available
at
 arXiv: math.GT/0008092

{\bf SkamB}
S. Kamada:
{\it Invariants of virtual braids and a remark on left stabilizations
and virtual exchange moves}, 
Kobe J. Math. {\bf 21} (2004), 33--49.





{\bf kauD}
L.~H.~Kauffman: 
{\it Virtual knot theory\/},
Europ.~J.~Combinatorics 
{\bf 20} (1999) 663--690. 






{\bf Kuper}
G. Kuperberg: 
{\it What is a virtual link?}, 
Algebr. Geom. Topol. {\bf 3} (2003), 587--591. 



{\bf Saw}
 J. Sawollek:   
{\it On Alexander-Conway polynomials for virtual knots and links\/}, 
preprint (1999), 
available
at
 arXiv: math.GT/9912173

{\bf SWb} D. S. Silver and S. G. Williams,  
{\it Alexander groups and virtual links\/}, 
preprint. 

{\bf SWA} 
D. Silver and S. Williams:
{\it Polynomial invariants of virtual links},
J. Knot Theory Ramifications {\bf 12} (2003), 987--1000. 


\bye

\end{thebibliography}
\bye